\documentclass{base}

\pdfoutput=1

\title{A 3-to-1 cactus graph: Details}
\author{Avi Levy}

\u{tkz-graph}
\u{xstring} %

\rc{\v}{\Vertex}
\nc{\e}{\Edge}

\nc{\el}[3]{
  \Edge[label = $#3$](#1)(#2)
}

\nc{\graphbdry}{
  \v[x = 0, y = 2]{2}
  \v[x = 2, y = 2]{3}
  \v[x = 6, y = 2]{5}
  \v[x = -1, y = 1]{6}
  \v[x = 3, y = 1]{7}
  \v[x = 5, y = 1]{8}
  \v[x = -3, y = 0]{9}
  \v[x = -2, y = 0]{10}
  \v[x = 0, y = -1]{13}
  \v[x = 2, y = -1]{14}
  \v[x = 5, y = -1]{15}
  \v[x = 5, y = -2]{18}
}

\nc{\graphint}{
  \v[x = -2, y = 2]{1}
  \v[x = 4, y = 2]{4}
  \v[x = 1, y = 0]{11}
  \v[x = 6, y = 0]{12}
  \v[x = -2, y = -2]{16}
  \v[x = 4, y = -2]{17}
}

\nc{\graphedges}{
  \foreach \v in {2,6,9,10}{
    \e(1)(\v);
  }
  \node (2 6) at (-.5,1.5) {};
  \node (9 10) at (-2.5,0) {};
  \node (13 14) at (1,-1) {};

  \draw[l] (13 14) to (2 6) to (9 10) to [bend right=60] (13 14);

  \node (7 14) at (2.5, 0) {};
  \node (5 8) at (5.5,1.5) {};
  \node (14 15) at (3.5, -1) {};

  \draw[r] (13 14) to (7 14) to (5 8) to [bend left=60] (14 15) to [bend left=60] (13 14);
  
  \foreach \v in {9,10,13,14}{
    \e(16)(\v);
  }
  
  \foreach \v in {2,3,6,7,13,14}{
    \e(11)(\v);
  }
  
  \foreach \v in {5,7,8,14}{
    \e(4)(\v);
  }
    
  \foreach \v in {5,8,14,15}{
    \e(12)(\v);
  }

  \foreach \v in {13,14,15,18}{
    \e(17)(\v);
  }
}

\rc{\d}[1]{
\IfStrEqCase{#1}{
  {0}{\EA}
  {1}{\NOEA}
  {2}{\NO}
  {3}{\NOWE}
  {4}{\WE}
  {5}{\SOWE}
  {6}{\SO}
  {7}{\SOEA}
}[]
}

\newcommand{\MRhref}[2]{\href{http://www.ams.org/mathscinet-getitem?mr=#1}{MR#2}}

\makeatletter
\def\@strippedMR{}
\def\@scanforMR#1#2#3\endscan{%
  \ifx#1M\ifx#2R\def\@strippedMR{#3}%
  \else\def\@strippedMR{#1#2#3}%
  \fi\fi}
\def\@rst#1 #2other{#1}
\newcommand\MR[1]{\relax\ifhmode\unskip\spacefactor3000 \space\fi
  \@scanforMR#1\endscan
  \MRhref{\expandafter\@rst\@strippedMR other}{\@strippedMR}}
\newcommand\MRs[1]{\relax\ifhmode\unskip\spacefactor3000 \space\fi
  \@scanforMR#1\endscan
  \MRhref{\@strippedMR}{\@strippedMR}}
\makeatother

\begin{document}
\maketitle
\b{abstract}{
  We briefly present a new construction of a 3-to-1 unrecoverable electrical network. See \cite{CIM}, \cite{French} and \cite{Klumb} for more background, and see \cite{Bireg} for a related construction. This is a taste of a comprehensive forthcoming paper on unrecoverable electrical networks.
}

\tikzset{
  base/.style = {
    thick, ->, shorten > = 1ex
  }
}
\tikzset{
  l/.style = {
    base, purple, bend right
  }
}
\tikzset{
  r/.style = {
    base, blue, bend left
  }
}
\tikzset{
  a/.style = {
    orange
  }
}

\fig{\centering
\b{tikzpicture}{[x = 2cm, y = 2cm, > = latex]
  \SetGraphUnit{2}

  \SetVertexSimple[FillColor = white]
  \graphint
  
  \SetVertexSimple[FillColor = black]
  \graphbdry

  \graphedges  
}
}{
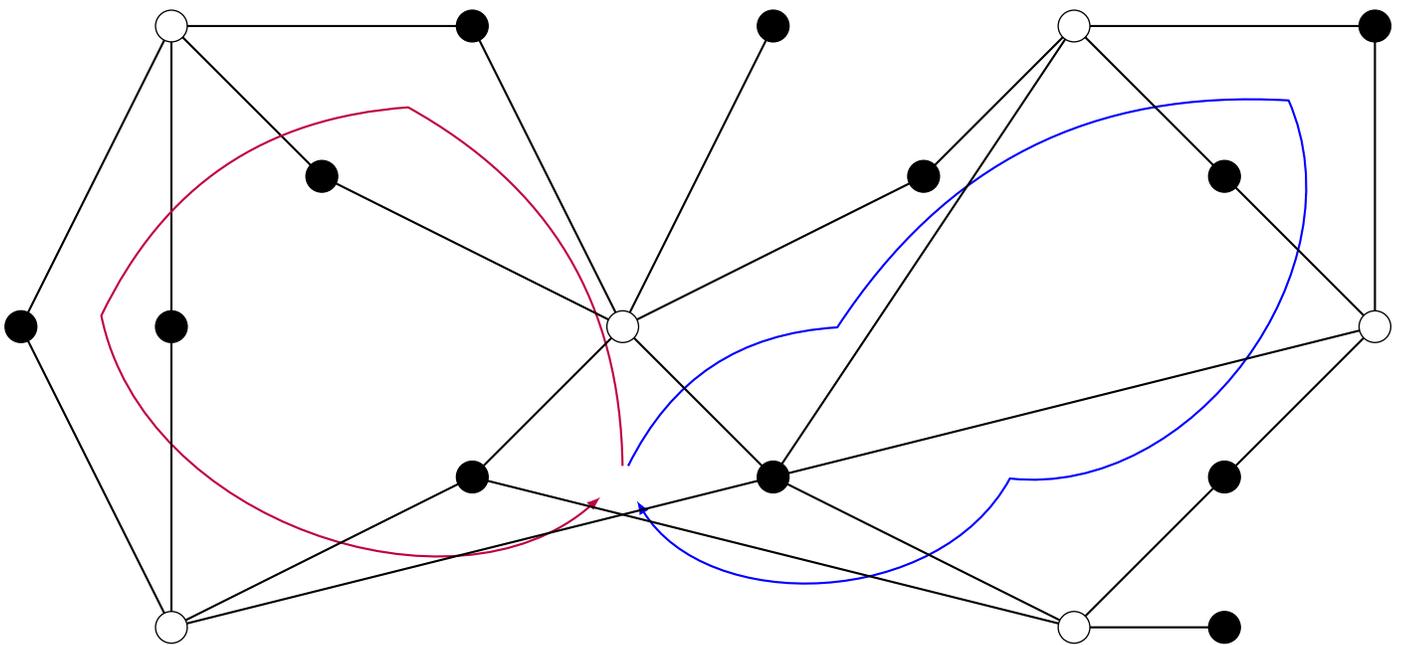
  Overview of a two leaf cactus: $\bullet$ = boundary vertex, $\circ$ = interior vertex.
}

The overview emphasizes the underlying loop structure. The left loop consists of quad$^3$, and the right consists of switch quad$^2$ switch. Note that the overview omits some ``non-structural'' auxiliary edges which we introduce on the next page.

Notice that the loops emanate from the central 6-star, which is actually a (quad \# switch). We also refer to it as a {\bf multiplexor}, following \cite{Klumb}.

\newpage

\fig{\centering
\b{tikzpicture}{[x = 2cm, y = 2cm, > = latex]
  \SetGraphUnit{2}

  \SetVertexSimple[FillColor = white]
  \graphint
  
  \SetVertexSimple[FillColor = black]
  \graphbdry

  \foreach \v in {2,6,9,10}{
    \e(1)(\v);
  }
  
  \foreach \v in {9,10,13,14}{
    \e(16)(\v);
  }

  \foreach \v in {2,3,6,7,13,14}{
    \e(11)(\v);
  }

  \foreach \v in {5,7,8,14}{
    \e(4)(\v);
  }
  
  \foreach \v in {5,8,14,15}{
    \e(12)(\v);
  }

  \foreach \v in {13,14,15,18}{
    \e(17)(\v);
  }

  \draw[a]
    (7) -- (6) -- (3) (6) -- (13)
    (2) -- (14) -- (3) (14) -- (18);
}
}{
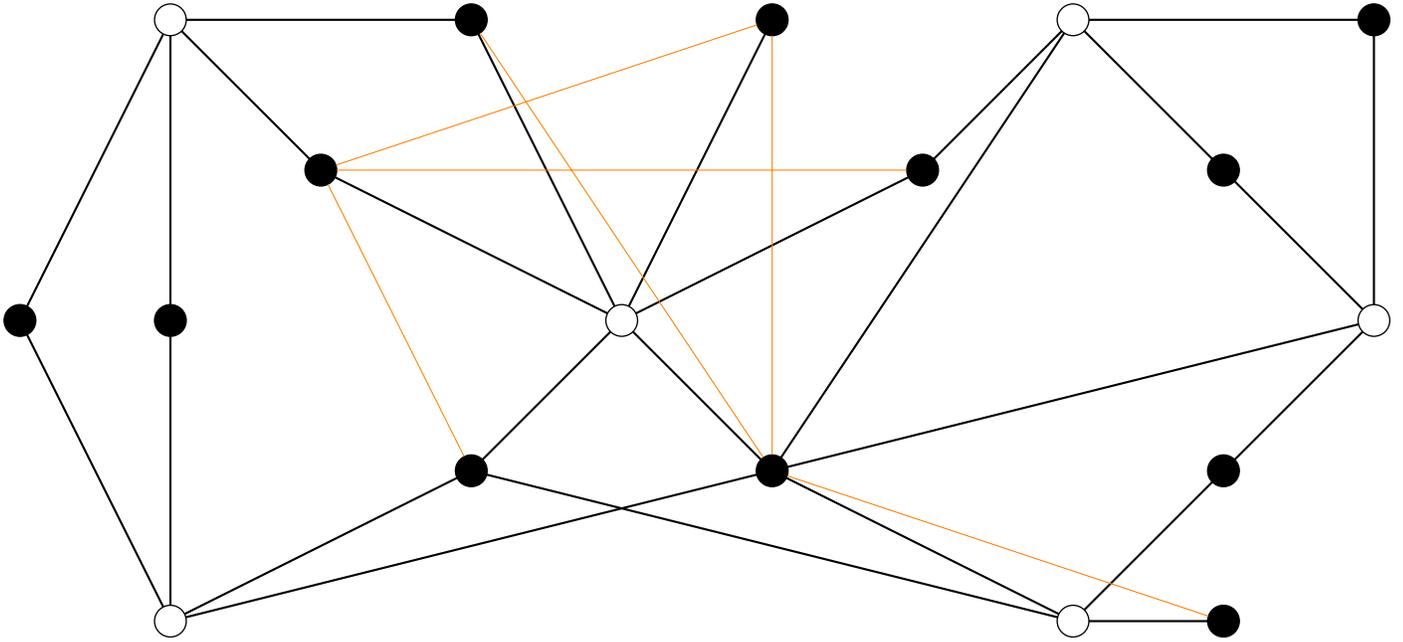
  Auxiliary edges which prevent recoverability
}

Here $\kappa(G)$ denotes the arity of $G$ (cardinality of largest fiber).
\b{claim}{[Upper bound]
  $\kappa(G)\leq 3$
}
\b{proof}{
Let $x$ denote the sum of the conductivities entering $v_{13}v_{14}$ (see Figure \ref{vertex number}). That is, $x$ is the sum of the red and blue arrow heads.

Suppose $x$ is known. By the claim following Figure \ref{quad switch}, all edges in the multiplexor are determined. Propagation along the red and blue loops determines the remaining edges. Hence
\[
\kappa(G)\leq \mbox{\# of choices for $x$}.
\]
Since there are 2 loops, $x$ satisfies a cubic and thus the graph is at most 3-to-1.
}
\b{claim}{[Lower bound]
$\kappa(G)\geq 3$
}
\b{proof}{
  We construct three networks on $G$ with the same response. The networks correspond to $x=2,3,4$. See tables on next page and subsequent diagrams for the rest of the construction.
}

\newpage

\fig{
\begin{tikzcd}[row sep = small, column sep = small]
  x\ar{r} & 7-x\ar{r} & \f{1}{7-x}\ar{r} & 2-\f{1}{7-x} \ar{r} & \f{6}{2-\f{1}{7-x}}\ar{r} & 4-\f{6}{2-\f{1}{7-x}}\ar{r} & \f{1}{4-\f{6}{2-\f{1}{7-x}}} & (=1-\f{3/2}{x-5})\\
  2 & 5 & \f{1}{5} & \f{9}{5} & \f{10}{3} & \f{2}{3} & \f{3}{2}\\
  3 & 4 & \f{1}{4} & \f{7}{4} & \f{24}{7} & \f{4}{7} & \f{7}{4}\\
  4 & 3 & \f{1}{3} & \f{5}{3} & \f{18}{5} & \f{2}{5} & \f{5}{2}\\
\end{tikzcd}
}{
  Arm propagation for \textcolor{purple}{left loop} (quad$^3$)
}
\fig{
\begin{tikzcd}[row sep = small, column sep = small]
  x\ar{r} & 7-x\ar{r} & 7-x\ar{r} & x \ar{r} & \f{1}{x}\ar{r} & 1-\f{1}{x}\ar{r} & \f{3/2}{1-\f{1}{x}}\ar{r} & \f{7}{2}-\f{3/2}{1-\f{1}{x}}\ar{r} & \f{7}{2}-\f{3/2}{1-\f{1}{x}} & (=2-\f{3/2}{x-1})\\
  2 & 5 & 5 & 2 & \f{1}{2} & \f{1}{2} & 3 & \f{1}{2} & \f{1}{2}\\
  3 & 4 & 4 & 3 & \f{1}{3} & \f{2}{3} & \f{9}{4} & \f{5}{4} & \f{5}{4}\\
  4 & 3 & 3 & 4 & \f{1}{4} & \f{3}{4} & 2 & \f{3}{2} & \f{3}{2}\\
\end{tikzcd}
}{
  Arm propagation for \textcolor{blue}{right loop} (switch quad$^2$ switch)
}
Observe that all three assignments are valid due to ``loop conservation'':
\[
\b{bmatrix}{
  \dfrac{3}{2}\\[1em]
  \dfrac{7}{4}\\[1em]
  \dfrac{5}{2}\\
}+\b{bmatrix}{
  \dfrac{1}{2}\\[1em]
  \dfrac{5}{4}\\[1em]
  \dfrac{3}{2}\\
}=
\b{bmatrix}{
  2\\
  3\\
  4\\
}
\]

\newpage

In a sense the construction is complete, for we have found a valid loop assignment. However, it is instructive to illustrate the process of populating a graph given a loop assignment. Our present graph can be decomposed into quads, switches, and a (quad \# switch), so it suffices to populate these three subgraphs. In the following diagrams, the expressions in the central interior node is a multiplier, to be applied to the surrounding edge weights.
\fig{
  \begin{minipage}{.5\textwidth}
    \centering
    \b{tikzpicture}{[x = 2cm, y = 2cm, > = latex]
      \SetGraphUnit{2}
      \SetVertexMath
      \SetVertexNormal[FillColor = white]

      \v[L = \f{1}{s}+2+t, x = 0, y = 0]{1}

      \SetVertexSimple[FillColor = black]
      \v[x = 1, y = 1]{2}
      \v[x = -1, y = 1]{3}
      \v[x = 1, y = -1]{4}
      \v[x = -1, y = -1]{5}

      \el{1}{2}{1}
      \el{1}{3}{s}
      \el{1}{4}{st}
      \el{1}{5}{s}

      \node (s) at (-1,0) {$s$};
      \node (t) at (1,0) {$t$};

      \draw (s) edge[l] (t);
    }
  \end{minipage}
  \begin{minipage}{.5\textwidth}
    \centering
    \b{tikzpicture}{[x = 2cm, y = 2cm, > = latex]
      \SetGraphUnit{2}
      \SetVertexMath
      \SetVertexNormal[FillColor = white]
        \v[L = s+2+\f{t}{s}, x = 0, y = 0]{1}

      \SetVertexSimple[FillColor = black]
      \v[x = -1, y = 1]{2}
      \v[x = 1, y = 1]{3}
      \v[x = -1, y = -1]{4}
      \v[x = 1, y = -1]{5}

      \el{1}{2}{1}
      \el{1}{3}{\f{t}{s}}
      \el{1}{4}{1}
      \el{1}{5}{s}
      
      \draw[a, bend left] (2) edge (5);

      \node (x) at (0,-1) {$s$};
      \node (y) at (1,0) {$t$};

      \draw (x) edge[r] (y);
    }
  \end{minipage}
}{
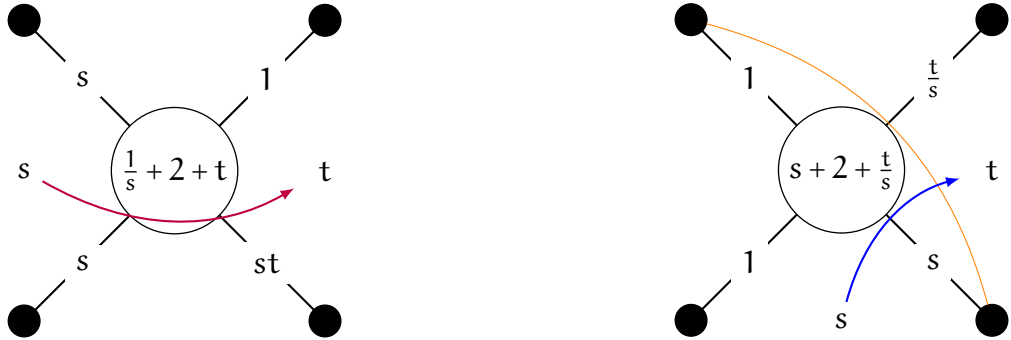
  Populating a quad and switch, respectively.
}

\fig{\centering
\b{tikzpicture}{[x = 2cm, y = 2cm, > = latex]
  \SetGraphUnit{2}

  \node (x) at (0,-1) {$s$};
  \node (z) at (1.5,0) {$t_2$};
  \node (y) at (-1.5,1.5) {$t_1$};

  \draw (x) edge[r] (z);
  \draw (x) edge[l] (y);

  \SetVertexNormal[FillColor = black]
  \SetVertexNoLabel
  \tikzset{
    VertexStyle/.append style = {minimum size = \VertexSmallMinSize}
  }
  \v[x = -1, y = 2]{2}
  \v[x = 1, y = 2]{3}
  \v[x = -2, y = 1]{4}
  \v[x = 2, y = 1]{5}
  \v[x = -1, y = -1]{6}
  \v[x = 1, y = -1]{7}

  \draw[a] 
    (5) -- (4) -- (6)
    (2) -- (7) -- (3) -- (4);
  
  \SetVertexMath
  \SetVertexNormal[FillColor = white]
  \SetVertexLabel

  \v[L = s+t_1+3+\f{t_2}{s}, x = 0, y = 0]{1}

  \el{1}{2}{1}
  \el{1}{3}{1}
  \el{1}{4}{t_1}
  \el{1}{5}{\f{t_2}{s}}
  \el{1}{6}{1}
  \el{1}{7}{s}

}}{
  \label{quad switch} Populating a (quad \# switch)
}

\b{claim}{
  The edge entering a (quad \# switch) determines all edges.
}
\b{proof}{
  Play the following game: an orange edge is removed if its endpoints can be connected with white edges. Observe that all orange edges are removed.
}
\newpage

\fig{\centering
\b{tikzpicture}{[x = 2cm, y = 2cm, > = latex]
  \SetGraphUnit{2}

  \SetVertexSimple[FillColor = white]
  \graphint
  
  \SetVertexSimple[FillColor = black]
  \graphbdry

  \el{1}{2}{\f{53}{5}}
  \el{1}{6}{\f{53}{5}}
  \el{1}{9}{\f{106}{3}}
  \el{1}{10}{\f{53}{9}}
  
  \el{16}{9}{\f{10}{3}}
  \el{16}{10}{\f{10}{3}}
  \el{16}{13}{5}
  \el{16}{14}{5}

  \el{11}{2}{\f{46}{5}}
  \el{11}{3}{\f{46}{5}}
  \el{11}{6}{\f{46}{25}}
  \el{11}{7}{\f{46}{5}}
  \el{11}{13}{\f{46}{5}}
  \el{11}{14}{46}

  \el{4}{5}{3}
  \el{4}{7}{6}
  \el{4}{8}{3}
  \el{4}{14}{6}
  
  \el{12}{5}{\f{7}{2}}
  \el{12}{8}{\f{7}{2}}
  \el{12}{14}{7}
  \el{12}{15}{\f{21}{2}}

  \el{17}{13}{\f{7}{2}}
  \el{17}{14}{\f{7}{4}}
  \el{17}{15}{\f{7}{2}}
  \el{17}{18}{\f{7}{2}}
}
}{
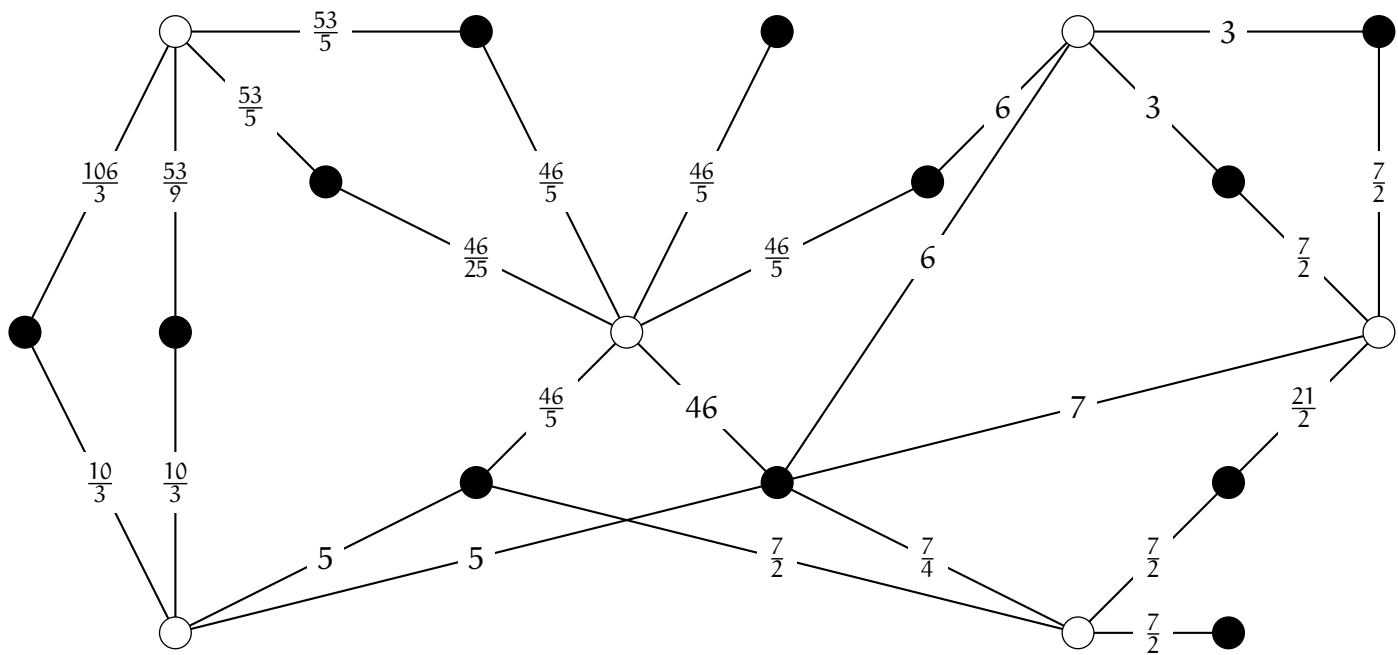
  Populated with $x=2$
}
\fig{\centering
\b{tikzpicture}{[x = 2cm, y = 2cm, > = latex]
  \SetGraphUnit{2}

  \SetVertexSimple[FillColor = white]
  \graphint
  
  \SetVertexSimple[FillColor = black]
  \graphbdry

  \el{1}{2}{\f{21}{2}}
  \el{1}{6}{\f{21}{2}}
  \el{1}{9}{36}
  \el{1}{10}{6}
  
  \el{16}{9}{\f{22}{7}}
  \el{16}{10}{\f{22}{7}}
  \el{16}{13}{\f{11}{2}}
  \el{16}{14}{\f{11}{2}}

  \el{11}{2}{\f{33}{4}}
  \el{11}{3}{\f{33}{4}}
  \el{11}{6}{\f{33}{16}}
  \el{11}{7}{\f{33}{4}}
  \el{11}{13}{\f{33}{4}}
  \el{11}{14}{33}

  \el{4}{5}{\f{8}{3}}
  \el{4}{7}{8}
  \el{4}{8}{\f{8}{3}}
  \el{4}{14}{8}
  
  \el{12}{5}{\f{23}{6}}
  \el{12}{8}{\f{23}{6}}
  \el{12}{14}{\f{23}{4}}
  \el{12}{15}{\f{69}{8}}

  \el{17}{13}{\f{17}{4}}
  \el{17}{14}{\f{85}{18}}
  \el{17}{15}{\f{17}{4}}
  \el{17}{18}{\f{17}{4}}
}
}{
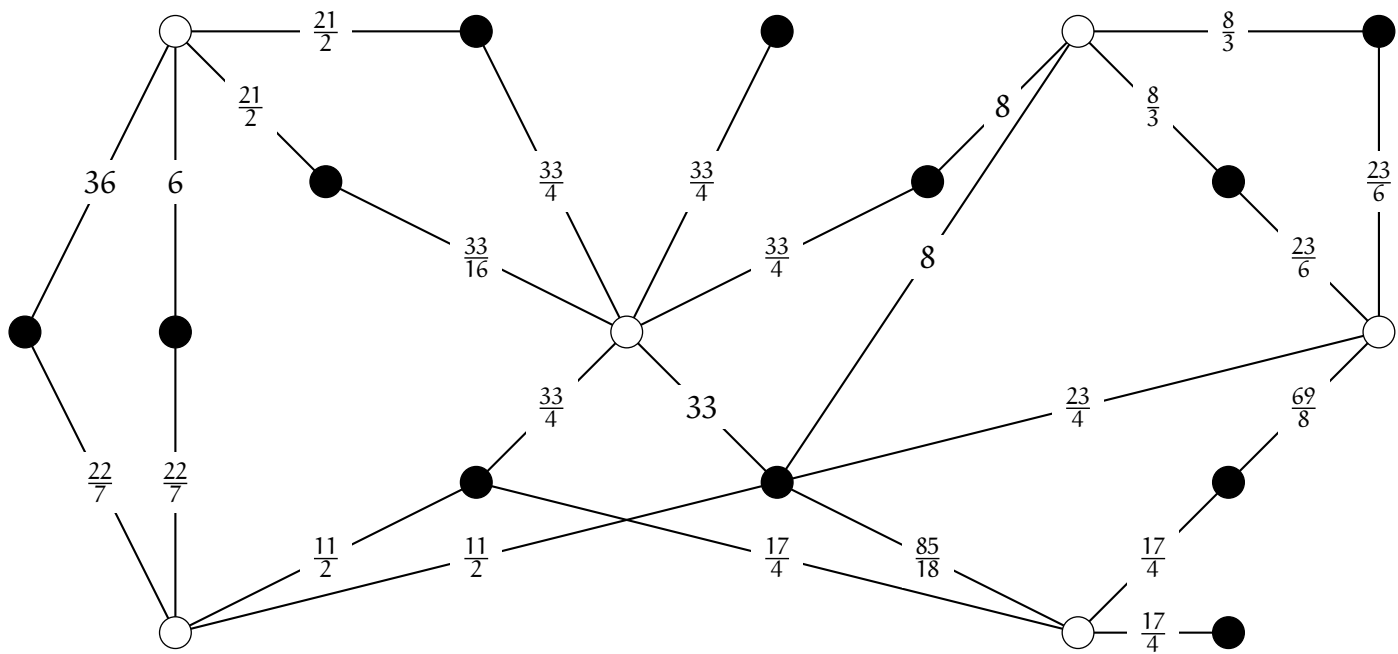
  Populated with $x=3$
}
\fig{\centering
\b{tikzpicture}{[x = 2cm, y = 2cm, > = latex]
  \SetGraphUnit{2}

  \SetVertexSimple[FillColor = white]
  \graphint
  
  \SetVertexSimple[FillColor = black]
  \graphbdry

  \el{1}{2}{\f{31}{3}}
  \el{1}{6}{\f{31}{3}}
  \el{1}{9}{\f{186}{5}}
  \el{1}{10}{\f{31}{5}}
  
  \el{16}{9}{\f{14}{5}}
  \el{16}{10}{\f{14}{5}}
  \el{16}{13}{7}
  \el{16}{14}{7}

  \el{11}{2}{\f{22}{3}}
  \el{11}{3}{\f{22}{3}}
  \el{11}{6}{\f{22}{9}}
  \el{11}{7}{\f{22}{3}}
  \el{11}{13}{\f{22}{3}}
  \el{11}{14}{22}

  \el{4}{5}{\f{5}{2}}
  \el{4}{7}{10}
  \el{4}{8}{\f{5}{2}}
  \el{4}{14}{10}
  
  \el{12}{5}{4}
  \el{12}{8}{4}
  \el{12}{14}{\f{16}{3}}
  \el{12}{15}{8}

  \el{17}{13}{\f{9}{2}}
  \el{17}{14}{\f{27}{4}}
  \el{17}{15}{\f{9}{2}}
  \el{17}{18}{\f{9}{2}}
}
}{
  Populated with $x=4$
}

Note that we didn't label the auxiliary edges, because it is clear that auxiliary edges can be chosen arbitrarily to ensure proper conductivities.

\fig{\centering
\b{tikzpicture}{[x = 2cm, y = 2cm, > = latex]
  \SetGraphUnit{2}

  \SetVertexLabel

  \graphint
  \graphbdry

  \graphedges
}
}{
  \label{vertex number} Standard vertex labelling
}

\newpage

\bibliographystyle{hmralpha}

\end{document}